\documentclass[a4paper,10pt]{article}

\usepackage[english]{babel}
\usepackage[T1]{fontenc}
\usepackage{enumerate}
\usepackage{amsmath}
\usepackage{amsthm}
\usepackage{amssymb}

\newcommand{\Qp}{\mathbb{Q}_p}

\newcommand{\Zp}{\mathbb{Z}_p}

\newcommand{\Q}{\mathbb{Q}}
\newcommand{\R}{\mathbb{R}}

\newcommand{\C}{\mathbb{C}}
\newcommand{\Z}{\mathbb{Z}}
\newcommand{\N}{\mathbb{N}}

\newcommand{\Square}[2]{\,\square_{#1}^{#2}\,}

\newtheorem{theorem}{Theorem}[section]

\newtheorem{prop}[theorem]{Proposition}
\newtheorem{lemma}[theorem]{Lemma}
\newtheorem{cor}[theorem]{Corollary}
\newtheorem{definition}[theorem]{Definition}
\newtheorem{Fact}{Fact}
\newtheorem{claim}[theorem]{Claim}

\theoremstyle{remark}
\newenvironment{Remark}{\begin{trivlist}\item[\hskip \labelsep {\bfseries Remark.}]}{\end{trivlist}}

\newenvironment{notation}{\begin{trivlist}\item[\hskip \labelsep {\bfseries Notations.}]}{\end{trivlist}}

\title{Expansions of the $p$-adic numbers that interprets the ring of integers.}
\author{Nathana\"el Mariaule\footnote{During the preparation of this paper the author was supported by the Fonds de la Recherche Scientifique - FNRS}} 
\date{}
\begin{document}
\maketitle
\begin{abstract} Let $\widetilde{\Qp}$ be the field of $p$-adic numbers in the language of rings.
In this paper we consider the theory of $\widetilde{\Qp}$ expanded by two predicates interpreted by multiplicative subgroups $\alpha^\Z$ and $\beta^\Z$ where $\alpha, \beta\in\N$ are multiplicatively independent. We show that the theory of this structure interprets Peano arithmetic if $\alpha$ and $\beta$ have positive $p$-adic valuation. If either $\alpha$ or $\beta$ has zero valuation we show that the theory of $(\widetilde{\Qp}, \alpha^\Z, \beta^\Z)$ does not interpret Peano arithmetic. In that case we also prove that the theory is decidable iff the theory of $(\widetilde{\Qp}, \alpha^\Z\cdot \beta^\Z)$ is decidable.
\end{abstract}

Questions about expansions of structures by powers of an integer have been around for a long time. In the 60's, B\"uchi proved that the theory of $(\Z,+,0, <, 2^\Z)$ is decidable. In a different spirit, L. van den Dries \cite{vdDries5} axiomatised the theory of the field of real numbers with a predicate for the powers of 2. More recently a $p$-adic equivalent for this latter result has been proved \cite{Mariaule2}. Having a good grasp of the expansion by one group, it is quite natural to look at the expansion by any collection of such groups. It turns out that this structure is much more complicated.
\par P. Hieronymi \cite{Hieronymi} proved that the theory of $(\R, +,\cdot, 2^\Z, 3^\Z)$ defines $\Z$ and therefore is undecidable. For the integers, it is not known whether the theory of $(\Z,+,2^\Z,3^\Z)$ is decidable or not. In this paper, we will discuss the question of decidability of $(\Qp,+,\cdot, \alpha^\Z, \beta^\Z)$ depending on $\alpha,\beta\in \N$ and $p$ prime number.
\par Let us remark that the group $\alpha^\Z$ has a different topological nature in $\Qp$ according to the valuation of $\alpha$. If $\alpha$ has positive $p$-adic valuation then $\alpha^\Z$ is a discrete group (isomorphic to a subgroup of the value group via the valuation). If $\alpha$ has zero valuation then it is dense in a finite union of multiplicative cosets of $1+p^k\Zp$ (where $k$ is the valuation of $\alpha-1$). We end up with three different cases: (1) if both $\alpha$ and $\beta$ have positive valuations. This is done in section \ref{discrete groups}. In that case $\alpha^\Z$ is in definable bijection with $(\alpha/\beta)^\Z$ and we get undecidability iff this latter group is dense in an open neighbourhood of $1$. Case (2): if $\alpha$ has positive valuation and $\beta$ has zero valuation. In that case an axiomatisation of the theory is given in \cite{Mariaule2}. The important ingredients of this axiomatisation are: first the axiomatisation of the theory of valued group induced on $\beta^\Z$, second the so-called Mann property of the group $\alpha^\Z\cdot \beta^\Z$ and smallness (see section \ref{dense groups} for the definitions), third the density of $\beta^\Z$ in a definable open neighbourhood of $1$  and last a definable bijection between $\alpha^\Z$ and a definable subgroup of the value group.
 Finally, case (3): if both $\alpha$ and $\beta$ have zero valuations. Here we adapt the strategy of case (2). First in section \ref{pair of valued groups} we look at a structure induced on the group $\alpha^\Z\cdot \beta^\Z$ i.e., we study the pair of groups $(\alpha^\Z\cdot \beta^\Z, \alpha^\Z)$ in a language of valued groups. Then we use it in section \ref{dense groups} to give an axiomatisation of the theory of $(\Qp,+,\cdot, \alpha^\Z, \beta^\Z)$. Again it is crucial that the group $\alpha^\Z \beta^\Z$ has the Mann property, is small and that both $\alpha^\Z$ and $\beta^\Z$ are dense in a definable open neighbourhood of $1$. In section \ref{NIP section} with the back-and-forth system used in the proof of the axiomatisation we give a description of definable sets. Then we show that the theory of $(\Qp,+,\cdot, \alpha^\Z, \beta^\Z)$ is NIP and therefore does not interprets Peano arithmetic if either $\alpha$ or $\beta$ has zero valuation. 

\begin{notation} $A^\times$ will denote the set of units in a ring. We denote the $p$-adic valuation by $v_p$. We always consider $\Qp$ with the language $\mathcal{L}_{Mac}=(+,-,\cdot, 0, 1, (P_n)_{n\in \N})$ where $P_n$ is interpreted by the set of $n$th powers. If $K$ is a valued field $K^h$ will denote its henselisation.
\end{notation}

\section{Expansion by two discrete groups}\label{discrete groups}
In this section we consider the case where the two subgroups are generated by elements $\alpha, \beta$ of positive $p$-adic valuation. In that case if $\alpha$ and $\beta$ are multiplicatively independent we obtain a definable bijection between $\alpha^\Z$ and a dense set. Using this and the structure of valued fields we obtain that the ring of integer is interpretable in our theory. Let us remark that Hieronymi proved in the real case that a definable bijection between any definable discrete infinite set and definable dense set implies that $\Z$ is definable \cite{Hieronymi}.

\begin{theorem}
Let $\alpha, \beta\in \N$ nonzero with $v_p(\alpha),v_p(\beta)>0$. Then, $Th(\Qp, +,\cdot, \alpha^\Z,\beta^\Z)$ is decidable iff $\alpha^\Z=\beta^\Z\not=\{1\}$.
\end{theorem}
\begin{proof}
First if $(\alpha)^\Z\cap(\beta)^\Z=\gamma^\Z$ for some $\gamma\not=1$ then $Th(\Qp, +,\cdot, \alpha^\Z,\beta^\Z)$ is interdefinable with $Th(\Qp, +,\cdot, \gamma^\Z)$ (for $\gamma^\Z$ is a subgroup of finite index in $\alpha^\Z$ and  $\beta^\Z$. The theory of this latter structure is decidable by \cite{Mariaule2}. Indeed Theorem 2.2 in that paper gives an axiomatisation of the theory. This axiomatisation is obviously recursively enumerable and therefore the theory is decidable.
\par Replacing $\alpha$ and $\beta$ by one of their power, we may assume that $v_p(\alpha)=v_p(\beta)$. Then $\gamma:= \alpha/\beta\in \Zp\setminus p\Zp$. Let us remark that $\gamma$ cannot be a root of unity by hypothesis on $\alpha,\beta$. Therefore $\gamma^\Z$ is not discrete. Again we replace $\alpha$ and $\beta$ by one of their power if necessary so that $\gamma\in 1+p\Zp$. We remark that we have definable isomorphisms between $\alpha^\Z$, $\beta^\Z$ and $\gamma^\Z$. For let $\tau:\alpha^\Z\rightarrow \beta^\Z$ which sends $\alpha^n$ to the unique element of $\beta^\Z$ with valuation $v_p(\alpha^n)$ (that is $\beta^n$) and $\sigma: \alpha^\Z\rightarrow \gamma^\Z: \alpha^n\longmapsto \alpha^n/\tau(\alpha^n)=\gamma^n$. Next we note that $v_p(\gamma^n-1)=v_p(log_p(\gamma^n))=v_p(n)+v_p(\gamma-1)$. We claim that the map $\alpha^n\rightarrow \alpha^{v_p(n)}$ from $\alpha^\N$ to itself is definable. Indeed we have that for all $n\in \N$ there are unique $\alpha^m\in \alpha^\N$ and $0\leq k<v_p(\alpha)$ such that $v_p(n)=v_p(\gamma^n-1)-v_p(\gamma-1)=v(\alpha^m)+k$. This latter condition is definable and therefore so is the map $\alpha^n\longmapsto (\alpha^m)^{v_p(\alpha)}\cdot\alpha^{k}=\alpha^{v_p(n)}$.
\par This proves that the structure $(\N,+,v_p,<)$ is definable in our theory where $v_p: \N\setminus\{0\}\rightarrow \N$. It remains to prove that the theory of this structure is undecidable. We remark that the exponentiation is definable in this structure:
$$p^n=\min\{k\in \N: n=v_p(k)\}.$$
Also the unary function $V_p(n)$ sending $n$ to the highest power of $p$ dividing $n$ is definable (it is $n\longmapsto k\in p^\N$ with $v_p(k)=v_p(n)$). Therefore the structure $(\N, +, V_p, p^x)$ is definable. The theory of this structure interprets the ring of integers by a result of  Elgot-Rabin \cite{Elgot-Rabin} and therefore is undecidable.
\end{proof}

\section{Expansion by two dense groups}
In this section we treat the case of $(\Qp, \alpha^\Z, \beta^\Z)$ where $v_p(\alpha)=v_p(\beta)=0$. To start with, we will assume that $\alpha^\Z, \beta^\Z\subset 1+p\Zp$ and $v_p(\alpha-1)=v_p(\beta-1)$. Note that the theory of the structure $(\Qp, \alpha^\Z)$ is axiomatised in \cite{Mariaule2}. The axiomatisation relies on the following observation: let $G$ be a multiplicative subgroup of $1+p^k\Zp$ ($k$ minimal for this property) then the $p$-adic valuation induces a structure of valued group on $G$. For let us recall that the $p$-adic logarithmic map $\log_p$ induces an isomorphism between $(1+p^k\Zp,\cdot)$ and $(p^{k}\Zp,+)$. So $\log_p(G)/p^k$ is a subgroup of the valued group $(\Zp,+,v_p)$.  We also have that $v_p(\log_p(1+px))=v_p(px)$ for all $x\in \Zp$.	Therefore $V_p:G\rightarrow \N\cup\{\infty\}:g\longmapsto v_p(g-1)-k$ is a valuation on $G$ and $(G, \cdot, V_p)\cong (\log_p(G)/p^k, +, v_p)$ as valued groups. An important step in \cite{Mariaule2} is an axiomatisation and a quantifier result for the theory of the structure $(G, \cdot, V_p)$. In the first part of this section we adapt this step to our setting. That is let $G=\alpha^\Z \beta^\Z$ (note that this group is definable in our language). Now we have extra-structures definable on $G$: e.g., $(G, \alpha^\Z, V_p)$. Using the symmetry of the problem it will not be necessary to look at $(G, \alpha^\Z, \beta^\Z, V_p)$. In section \ref{pair of valued groups} we will prove a quantifier elimination result and give an axiomatisation for the theory of the pair of valued groups $(G, \alpha^\Z, V_p)$. Then in section \ref{dense groups} we will use these results to axiomatise the theory of $(\Qp, \alpha^\Z, \beta^\Z)$. Finally in the last subsection we prove that the theory of $(\Qp, \alpha^\Z, \beta^\Z)$ is NIP. In particular it does not interpret Peano arithmetic.

\subsection{Pairs of $p$-valued groups}\label{pair of valued groups}
Let $G$ be a subgroup of $(\Zp,+)$. Then the $p$-adic valuation induces on $G$ a structure of $p$-valued group.
\begin{definition}
Let $(G,+,0_G)$ be an abelian group and $V:G\rightarrow \Gamma\cup\{\infty\}$ where $\Gamma$ is a totally ordered set with discrete order, no largest element and $\infty$ is an element such that $\infty>\gamma$ for all $\gamma\in \Gamma$. We say that $(G,+, V)$ is a \emph{$p$-valued group} if for all $x,y\in G$ and for all $n\in \Z$,
\begin{itemize}
	\item $V(x)=\infty$ iff $x=0_G$;
	\item $V(nx)=V(x)+v_p(n)$;
	\item $V(x+y)\geq \min\{V(x),V(y)\}$;
\end{itemize}
where $v_p$ is the $p$-adic valuation, $nx=x+\cdots +x$ ($n$ times), $(-n)x=-(nx)$ for all $n>0$, $0x=0_G$ and if $x\in G$, $V(x)+k$ denotes the $k$th successor of $V(x)$ in $\Gamma\cup\{\infty\}$ (by convention the successor of $\infty$ is $\infty$). 
\end{definition}
It is clear that $(\Zp, +, v_p)$ and $(\alpha^\Z\cdot \beta^\Z, \cdot, V_p)$ are $p$-valued groups.
In this section we consider a pair of $p$-valued groups $(G,H,V)$ (i.e., $H$ is a subgroup of $G$ and the valuation on $H$ is the restricted valuation from $G$) such that
\begin{itemize}
	\item $[q]H:= [H:qH]= q$ for all prime $q$;
	\item $[q]G=q^2$ for all prime $q$;
	\item $G/H$ is torsion-free, infinite;
	\item $H$ is dense, codense in $G$.
\end{itemize}
Example of such groups are: $(\Z+x\Z, \Z)$ with the $p$-adic valuation where $x\in \Qp\setminus \Q$ or $(\alpha^\Z\beta^\Z, \alpha^\Z, V_p)$ where $\alpha, \beta\in 1+p\Zp$, multiplicatively independent and $V_p(x)=v_p(x-1)-\min\{v(\alpha-1),v(\beta-1)\}$ (in this section we will assume that both elements in the min are equal). The theory of the valued groups $(G, V), (H,V)$ has been axiomatised in \cite{Mariaule2}. Also a quantifier result is proved (\cite{Mariaule2} Theorem 1.2). Here we will prove that adding the density axiom and purity assumption to the theory of these group is sufficient to treat the case of pair of groups.
\par We define $T_{pV}^{pair}$ the theory whose models $(G, H,+,-,0,1, C, \equiv_n, (VG\cup\{\infty\}, <,S,0_\Gamma, \infty), V)$ satisfy:

\begin{enumerate}
	\item $(G,+,-,0)$ is an abelian group, $x\equiv_n y$ iff $\exists g\in G, x=y+ng$ and $[q]G=q^{2}$ for all $q$ prime;
	\item $H$ is a pure subgroup of $G$, $[q]H=q$ for all $q$ prime and $1\in H$;
	\item $(VG,<)$ is a discrete ordered set with first element $0_{\Gamma}$, any nonzero element has a predecessor and there is no last element, $\infty$ is an element such that $\gamma<\infty$ for all $\gamma\in VG$ and $S$ is the successor function ($S(\infty):=\infty$);
	\item $(G,V)$, $(H,V)$ are $p$-valued groups;
	\item $1, C$ are elements such that $V(1)=V(C)=0_{\Gamma}$, $i\cdot 1+j\cdot C \not\equiv_n i'\cdot 1+j'\cdot C$ for all $(i,j)\not=(i',j')\in \{0,\cdots, n-1\}^2$. Also $i\cdot 1$ and $i'\cdot 1$ are in distinct cosets of $nH$ for all $i\not= i'\in \{0,\cdots, n-1\}$;
	\item For all $x,y\in G$, if $V(x)=V(y)$, then there is a unique $0<i<p$ such that $V(x-iy)>V(x)$;
	\item $G$ is regularly dense i.e., for all $n$ $nG$ is dense in $\{x\in G\mid\ V(x)\geq v_p(n)\}$ (where $v_p(n)$ denotes the $v_p(n)$th successor of $0_{\Gamma}$ in $VG$) i.e. for all $n$
	$$\forall x\in G\ V(x)\geq v_p(n)\rightarrow\Big[\forall\gamma\geq v_p(n)\in VG\exists y\in nG\ V(x-y)\geq \gamma\Big];$$
	\item for all $n\in \N$, $nH$ is dense, codense in $nG$.
\end{enumerate}
\begin{Remark}From the axioms, it follows that $G$ is torsion-free (as it is a $p$-valued group) and that for all $x,y\in H$, $x\equiv_n y$ iff there is $h\in H$ such that $x=y+nh$. Also, both $(G,+, 0,1,V)$ and $(H,+,0,1,V)$ are $p$-valued groups, models of theories described in \cite{Mariaule2} section 2. It is clear that $(\alpha^\Z\cdot \beta^\Z, \alpha^\Z, \cdot, V_p)$ is a model of the theory where $1$ is interpreted by $\alpha$ and $C$ by $\beta$.
\end{Remark}

\begin{theorem}\label{pair of $p$-valued groups} The theory $T_{pV}^{pair}$ admits the elimination of quantifiers.
\end{theorem}
From this theorem it follows that
\begin{cor} $Th(\alpha^\Z\beta^\Z, \alpha^\Z,V_p)$ is axiomatised by $T_{pV}^{pair}\cup tp(\beta/\alpha^\Z)$.
\end{cor}
For it is sufficient to remark that $(\alpha^\Z\beta^Z,\alpha^\Z, V_p)$ is a prime model. Then by quantifier elimination, $T_{pV}^{pair}\cup tp(\beta/\alpha^\Z)$ is complete.

We will now prove Theorem \ref{pair of $p$-valued groups}:
\begin{proof}
Let $M^*=(G^*,H^*)$ and $M'^*=(G'^*, H'^*)$ be saturated models of our theory. Let $M=(G,H)$ and $M'=(G',H')$ be isomorphic substructures of $(G^*,H^*)$ and $(G'^*, H'^*)$ (of cardinality less than the saturation). We denote by $\iota$ the isomorphism. An isomorphism between substructures (which are torsion-free groups) extends uniquely to pure closure (the language contains congruence relations). So we may assume that all subgroup inclusions are pure. Let $x^*\in M^*\setminus M$. We will prove that the isomorphism extend to $M\langle x^*\rangle$, the pure closure of $M(x^*)$ in $M^*$ i.e., to $\{g^*\in G^*\mid\ ng^*=mx^*+g$ for some $m,n\in \N$ and $g\in G\}$.
\par Let $\Phi(\overline{a},\overline{b})$ be the set of formulas of the form $V(lx-a_i)+r\Square{}{}V(b_i)$ such that $V(lx^*-a_i)+r\Square{}{}V(b_i)$ holds where $l,r\in \Z$, $\overline{a},\overline{b}\subset M$ and $\Square{}{}$ holds for $<,>,\leq, \geq$ or $=$. Let $\Psi(\overline{a})$ be the set of formulas of the form $lx-a_i\equiv_n0$ such that $lx^*-a_i\equiv_n0$ holds where $l\in \Z, n\in \N$ and $\overline{a}\subset M$.
\begin{claim}\label{type pair of valued group} For all $\overline{a},\overline{b}\subset M$, for all $\varphi_1(x,\overline{a},\overline{b}),\cdots, \varphi_k(x,\overline{a},\overline{b})\in \Phi$ and for all $\psi_1(x,\overline{a}),\cdots, \psi_l(x,\overline{a})\in \Psi(\overline{a})$, there is $y^*\in M'^*$ such that $\bigwedge_i\varphi_i(y^*,\iota(\overline{a}),\iota(\overline{b}))\wedge\bigwedge_j \psi_j(y^*,\iota(\overline{a}))$ holds. Furthermore, we may assume that $y^*\in H'^*$ iff $x^* \in H^*$. 
\end{claim}
First by the properties of congruences, $\bigwedge_j \psi_j(x,\overline{a})$ is equivalent to $x\equiv_n a$ for some $n\in \N$, $a\in M$. Note that if $x^*\in H^*$, we may assume $a\in H$ (even in $\{0,\cdots, (n-1)\cdot 1)\}$ by axioms 5.) . By the axiom of regular density, $a+nG^*$ is dense in $B(a, v_p(n))$ (the ball of centre $a$ and radius $v_p(n)$). So $x^*$ realises the formula  $\bigwedge_i\varphi_i(x,\overline{a},\overline{b})\wedge V(x-a)>v_p(n)$. By \cite{Mariaule2} Lemma 1.4, there is $y^*_1\in G'^*$ such that  $\bigwedge_i\varphi_i(y^*_1,\iota(\overline{a}),\iota(\overline{b}))\wedge V(y^*_1-\iota(a))>v_p(n)$. By properties of the valuation, there is an open ball $B$ around $y^*_1$ such that any point in $B$ satisfies the same formula. As $B\subset B(\iota(a),v_p(n))$ and $\iota(a)+nG'^*$ is dense in $B(\iota(a),v_p(n))$, there is $y^*$ such that $y^*\equiv_n\iota(a)$ and $y^*\in B$ (so $\bigwedge_i\varphi_i(y^*,\iota(\overline{a}),\iota(\overline{b}))$). Furthermore, as $nH'^*$ is dense, codense in $nG'^*$ according to the situation of $x^*$, we can take $y^*\in H'^*$ or $y^*\in G'^*\setminus H'^*$. This completes the proof of the claim.

\par By the above claim and saturation, there is $y^* \in G'^*$ which realises all formulas in $\Phi(\overline{a},\overline{b})$ and $\Psi(\overline{a})$ for all $\overline{a},\overline{b}\subset G^*$. Also we can take $y^*\in H'^*$ iff $x^*\in H^*$. Then $\iota$ extends to an isomorphism of valued groups between $G\langle x^*\rangle$ and $G'\langle y^*\rangle$ by $x^*\longmapsto y^*$. It remains to prove that for all $x\in G\langle x^*\rangle$, $x\in H^*$ iff $\iota(x)\in H'^*$. First if $x^*\in \langle GH^*\rangle$ i.e., $nx^*= g+h^*$ for some $g\in G, h^*\in H^*$ and $n\in \N$, then as $G\langle x^*\rangle=G\langle nx^*-g\rangle$, we may assume that $x^*\in H^*$. Let $x\in G\langle x^*\rangle$ i.e., $nx=mx^*+c$ for some $n,m\in \Z$ and $c\in G$. If $x\in H^*$ then $c\in H$. So $x\in H\langle x^*\rangle$ (the pure closure of $H(x^*)$ in $H^*$). But by choice of $y^*$ we also have that the extension of $\iota$ induces isomorphism between $H\langle x^*\rangle$ and $H'\langle y^*\rangle$. So we are done. Now assume that $x^*\notin \langle GH^*\rangle$. Then for all $x\in H^*$, $nx=mx^*+c$ for some $n,m\in \Z$ and $c\in G$ iff $m=0$ and $c\in H$. In particular for all $x\in G\langle x^*\rangle \cap H^*$, $nx\in H$ for some $n$. As $H$ is a pure subgroup of $H^*$ $x\in H$. This concludes the proof of the theorem.

\end{proof}

\subsection{Expansion by a pair of groups}\label{dense groups}
In this section we give an axiomatisation of the theory of $(\Qp, \alpha^{\Z}, \beta^\Z)$ where $v_p(\alpha)=v_p(\beta)=0$. First we introduce some definitions involved in this axiomatisation.

\par Let $K$ be a field of characteristic zero and $G$ be a subgroup of $K^\times$. Let $a_1,\cdots, a_n\in \Q$ nonzero. We consider the equation
$$a_1x_1+\cdots+a_nx_n=1.$$
A solution $(g_1,\cdots, g_n)$ of this equation in $G$ is called nondegenerate if $\sum_{i\in I}a_ig^i\not=0$ for all $I\subset\{1,\cdots, n\}$ nonempty. We say that $G$ has the \emph{Mann property} if for any equation like above there is finitely many nondegenerate solutions in $G$. Examples of groups with Mann property are the roots of unity in $\C$ \cite{Mann} or any group of finite rank in a field of characteristic zero (see \cite{Evertse-Gyory} Theorem 6.3.1 for instance). In particular, any subgroup of $\Qp^\times$ of finite rank has the Mann property, for instance it is the case for $\alpha^\Z \beta^\Z$.
\par Let $G<K^\times$ be a group with the Mann property. Then the Mann axioms are axioms in the language of rings expanded by constant symbols $\gamma_g$ for the elements of $G$ and a unary predicate $A$ for $G$. Let $a_1,\cdots , a_n\in\Q^\times$. As $G$ has the Mann property, there is a collection of $n$-uples $\overline{g}_i=(g_{1i},\cdots, g_{ni})$  ($1\leq i\leq l$) in $G^n$ so that these $n$-uples are the nondegenerate solutions of the equation $a_1x_1+\cdots+a_nx^n=1$. The corresponding Mann axiom express that there is no extra nondegenerate solution in $A$ i.e.,
$$\forall\overline{y}\left[\left(\bigwedge_i A(y_i)\wedge\sum_{i=1}^n a_iy_i=1\wedge\bigwedge_{I\subset \{1,\cdots ,n\}}\sum_{i\in I} a_iy_i\not=0\right)\rightarrow \bigvee_{k=0}^l \overline{y}=\overline{\gamma}_{g_k}\right]. $$
The main consequence of Mann axioms that we will use is the following:
\begin{lemma}[Lemmas 5.12 and 5.13 in \cite{Gunaydin-vdD}]\label{G-vdD Lemma 5.12} Let $K$ be a field of characteristic zero, let $G$ be a subgroup of $K^\times$ and let $\Gamma$ be a subgroup of $G$ such that for all $a_1,\cdots a_n\in \Q^\times$ the equation $a_1x_1+\cdots +a_nx_n$ has the same nondegenerate solutions in $\Gamma$ as in $G$. Then, for all $g, g_1,\cdots, g_n\in G$
\begin{itemize}
	\item if $g$ is algebraic over $\Q(\Gamma)$ of degree $d$ then $g^d\in \Gamma$;
	\item if $g_1,\cdots, g_n$ are algebraically independent over $\Q(\Gamma)$ then they are multiplicatively independent over $\Gamma$.
\end{itemize} 
In particular, if $\Gamma$ is  a pure subgroup of $G$, then the extension $\Q(G)$ over $\Q(\Gamma)$ is purely transcendental.
\end{lemma}

\par Let $\mathcal{M}=(M,\cdots)$ be a $\mathcal{L}$-structure. Let $A\subset M$ and $\mathcal{L}_A$ be the expansion of $\mathcal{L}$ by a unary predicate that will be interpreted by $A$ in $M$. We denote by $f: X\stackrel{n}{\rightarrow} Y$ a map from $X$ to the subsets of $Y$ of size at most $n$. We say that $A$ is \emph{large} in $M$ if there is a $\mathcal{L}_A$-definable map $f:M^m\stackrel{n}{\rightarrow} M$ such that $f(A)=\bigcup_{x\in A^m} f(x)=M$. We say that $G$ is \emph{small} if it is not large. Note that as $n^\Z\cdot m^\Z$ is countable, it is small in $\Qp$. Let us also remark that smallness can be written as a scheme of first-order sentences in the language $\mathcal{L}_A$.

Let $\alpha,\beta\in \N$ multiplicatively independent with $v_p(\alpha)=v_p(\beta)=0$ and $v_p(\alpha-1)=v_p(\beta-1)>0$. We set $\mathcal{L}_{G,H}$ to be the language $\mathcal{L}_{Mac}$ expanded by two unary predicates $G,H$ interpreted in $\Qp$ by $\alpha^\Z$ and $\beta^\Z$. We will now axiomatise the theory of $(\Qp, \alpha^\Z, \beta^\Z)$. Let $T_{\alpha,\beta}$ be the theory whose models $(K, G_K,H_K)$ satisfy:
\begin{itemize}
	\item $(K,+,-, \cdot, 0,1)$ is a $p$-adically closed field;
	\item $G_K$, $H_K$ are multiplicative subgroup of $K^\times$;
	\item $((G_K\cdot H_K, G_K, 1,\alpha,\beta, \equiv_k (k\in \N)), (vK\cup\{\infty\}, <,S, 0_\Gamma, \infty), V)$ is elementary equivalent to $((\alpha^\Z\cdot \beta^\Z, \alpha^\Z, 1,\alpha,\beta, \equiv_k (k\in\N)),(\N\cup\{\infty\}, <,S,< 0, \infty), V_p)$ where $V: G_K\cdot H_K\rightarrow vK\cup\{\infty\}: g\longmapsto v_K(g-1)-1$ and $g\equiv_n g'$ iff there is $z\in G_KH_K$ such that $g=g'z^n$. The axiomatisation of this structure is given in section \ref{pair of valued groups};
	\item $((G_K\cdot H_K,H_K, 1,\beta,\alpha, \equiv_k (k\in \N)), (vK\cup\{\infty\}, <,S, 0_\Gamma, \infty), V)$ is elementary equivalent to $((\alpha^\Z\cdot \beta^\Z, \beta^\Z, 1,\beta,\alpha, \equiv_k (k\in\N)),(\N\cup\{\infty\}, <,S,< 0, \infty), V_p)$;
	\item $G_K\cap H_K=\{1\}$;
	\item $G_K\cdot H_K$ satisfy the Mann axioms for $\alpha^\Z\cdot \beta^\Z$.
	\item $G_K,H_K$ are dense in $1+p^{v_p(\alpha-1)}\mathcal{O}_K$;
	\item $G_K\cdot H_K$ is a small set.
\end{itemize}
\begin{Remark} Note that the $p$-adic valuation is interpretable in the language of rings as $v_p(x)\geq 0$ iff $1+px^2$ has a square root in $\Qp$ (if $p\not=2$) or iff $1+px^3$ has a $3$rd root in $\Q_2$. Therefore the above set of axioms is expressible in the language $\mathcal{L}_{G,H}$.
\end{Remark}

\begin{theorem}\label{completeness}
$T_{G,H}$ is complete.
\end{theorem}
\begin{proof}
Let $(K^*, G_{K^*}, H_{K^*})$ and $(L^*, G_{L^*}, H_{L^*})$ be two saturated models of the theory with same cardinality. Let $Sub(K^*)$ be the collection of $\mathcal{L}_{G,H}$-substructures $(K', G_{K'}, H_{K'})$ of $(K^*, G_{K^*}, H_{K^*})$ such that 
\begin{itemize}
	\item $K'$ is $p$-adically closed, $|K'|<|K^*|$;
	\item  $G_{K'}H_{K'}$ (resp. $G_{K'},H_{K'}$) is a pure subgroup of $ G_{K^*}H_{K^*}$ (resp. $G_{K^*}, H_{K^*}$);
	\item  $K'$ and $\Q(G_{K^*}H_{K^*})$ are free over $\Q(G_{K'}H_{K'})$.
\end{itemize}
We define similarly $Sub(L^*)$. Note that as $\Q(G_{K^*}H_{K^*})$ is a regular extension of $\Q(G_{K'}H_{K'})$ (by Lemma \ref{G-vdD Lemma 5.12}) and by freeness  $K'$ and $\Q(G_{K^*}H_{K^*})$ are linearly disjoint over $\Q(G_{K'}H_{K'})$. We prove that these two sets and $\mathcal{L}_{G,H}$-isomorphisms between their elements have the back-and-forth property. First let us remark that $Sub(K^*)$ is nonempty. For $(\Q(\alpha^\Z,\beta^\Z)^h,\alpha^\Z,\beta^\Z)\in Sub(K^*)$. The same holds for $Sub(L^*)$.
We fix $(K', G_{K'}, H_{K'})\in  Sub(K^*)$, $(L', G_{L'}, H_{L'})\in  Sub(L^*)$ and $\iota$ an isomorphism between these structures. Let $x^*  \in K^*\setminus K'$. We shall prove that $\iota$ extends to an isomorphism having $x^* $ in its domain. There are 4 possible cases:
\par (1) If $x^* \in G_{K^*}$: let $p_K$ be the type of $x^* $ over $K'$ in the language of rings and $q_K$ its image by $\iota$. Let us remark that as $K^*, K'$ are $p$-adicallly closed $p_K$ is determined by the formulas of the form $v(x-a)\Square{}{} v(b)$ with $a,b\in K'$. Let $\phi_1(x,a_1,b_1),\cdots ,\phi_k(x,a_k,b_k)$ be a finite collection of these formulas. Let $p_G(x)$ be the type of $x^* $ over $(G_{K'}H_{K'}, G_{K'})$ (in the language of pair of $p$-valued $\Z$-groups) and let $q_G(x)$ be its image by $\iota$. By Claim \ref{type pair of valued group} this type is determined by formulas of the form $V(x-g)\Square{}{} v_k(a)$ and $x\equiv_n g$ for some $g\in G_{K'}H_{K'}$, $a\in K'$ and $n\in\N$.

 Then by the density axiom and the proof of quantifier elimination for the family of $p$-valued $\Z$-groups (Theorem \ref{pair of $p$-valued groups}), one can find a realisation of $q_G(x)\cup \{\phi_1(x,\iota(a_1),\iota(b_1)),\cdots ,\phi_k(x,\iota(a_k),\iota(b_k))\}$. So by saturation there is $y^*$ realisation of $q_K\cup q_G$.
 Let us note that $K'(x^* )^h\cong L'(y^*)^h$ as valued fields where the isomorphism $\iota'$ is the extension of $\iota$ by $x^* \longmapsto y^*$. We have that  $K'(x^* )^h\cap G_{K^*}H_{K^*}= G_{K'}H_{K'}\langle x^*  \rangle$ and $L'(y^*)^h\cap G_{L^*}H_{L^*}= G_{L'}H_{L'}\langle y^* \rangle$. This follows from the following fact:
\begin{Fact} Let $t_1,\cdots, t_n\in K'$ algebraically independent over $G_{K'}H_{K'}$. Then
$$acl(t_1,\cdots, t_n, x^* )\cap G_{K^*}H_{K^*}= G_{K'}H_{K'}\langle x^* \rangle,$$
where $acl$ is the algebraic closure relation in the language of rings.
\end{Fact}
This fact is a consequence of Mann property (Lemma \ref{G-vdD Lemma 5.12}), see \cite{Belegradek-Zilber} Lemma 4.2 for a proof.

As $y^*$ is a realisation of $q_G$, $\iota'$ induces an isomorphism between  $G_{K'}H_{K'}\langle x^*  \rangle$ and $G_{L'}H_{L'}\langle y^* \rangle$ in the language of pair of $p$-valued groups. Therefore for all $x\in K'(x^* )^h$, $x\in G_{K^*}$ iff $\iota(x)\in G_{L^*}$. It remains to prove that $x\in H_{K^*}$ iff $\iota'(x)\in H_{L^*}$. For if $x\in H_{K^*}$  and $x\in G_{K'}H_{K'}\langle x^*  \rangle$, there is $g\in G_{K'}$, $h\in H_{K'}$, $n,m\in \N$ such that $x^n={(x^*)}^mgh$. Therefore $x^nh^{-1}={(x^*)}^mg$. As $x^nh^{-1}\in H_{K^*}$ and ${(x^*)}^mg\in G_{K^*}$ this implies that $x^nh^{-1}\in H_{K^*}\cap G_{K^*}=\{1\}$. Therefore, $x^n=h'$ and as $H'$ is a pure subgroup of $H^*$ and $H^*$ is torsion-free, $x\in H'$. So $\iota'(x)=\iota(x)\in H_{L^*}$.

This proves that $(K'(x^* )^h, G_{K'}\langle x^* \rangle, H_{K'})\in Sub(K^*)$, $(L'(x^* )^h, G_{L'}\langle y^*\rangle, H_{L'})\in Sub(L^*)$ and $\iota'$ is an $\mathcal{L}_{G,H}$-isomorphism between these structures. This concludes this case.
\par (2) If $x^* \in H_{K^*}$: same as case (1).
\par (3) If $x^*  \in K'(G_{K^*}, H_{K^*})^h$: then $x^* \in K'(g_1,\cdots, g_k, h_1,\cdots, h_l)^h$ where $g_i\in G_{K^*}$ and $h_j\in  H_{K^*}$. This case follows from cases (1) and (2) by induction on $k,l$.
\par (4) If $x^*  \notin K'(G_{K^*}, H_{K^*})^h$: by smallness we can realise in $L^*$ any cut over $L'(G_{K^*}, H_{K^*})^h$. In particular let $y^*$ be a realisation of the image of the cut of $x^* $ over $K'(G_{K^*}, H_{K^*})^h$ by $\iota$. Then $(K'(x^* )^h, G_{K'}, H_{K'})\in Sub(K^*)$ and $(L'(y^*)^h, G_{L'}, H_{L'})\in Sub(L^*)$. Furthermore $\iota$ extends to an isomorphism between $K'(x^* )^h$ and $L'(y^*)^h$ with $x^* \longmapsto y^*$ by linear disjointness. This completes the proof of the theorem.

\end{proof}

\begin{Remark} From the proof of this theorem and the proof of Theorem 2.4 in \cite{Mariaule2} one can deduce an axiomatisation of $Th(\Qp, p^\Z, \alpha^\Z,\beta^\Z)$. For let $\mathcal{L}_p=\mathcal{L}_{G,H}\cup \{ A, \lambda\}$, where $A$ is a unary predicate interpreted in $\Qp$ by $p^\Z$ and $\lambda$ a function symbol interpreted by $x\longmapsto a\in A$ such that $v_p(x)=v_p(a)$. Let $T_p$ be the extension of $T_{G,H}$ by the following axioms:
\begin{itemize}
	\item $A$ is a multiplicative subgroup of $K^\times$, $p\in A$;
	\item $v_K$ induces a group isomorphism between $v_K$ and $A(K)$;
	\item $\forall x v_K(\lambda(x))=v_K(x)$ and $\lambda: K^\times \rightarrow A$ is surjective;
	\item Mann axioms for the group $p^\Z \alpha^\Z \beta^\Z$.
\end{itemize}
Then $T_p$ is a complete theory. This follows from the proof of Theorem 2.4 in \cite{Mariaule2}  where in step 1.(b) we use steps (1)-(2) from the proof of Theorem \ref{completeness}.
\end{Remark}

\begin{cor}\label{effective MP} Let $\alpha, \beta\in \N$ with $v_p(\alpha)=v_p(\beta)=0$. Then 
$Th(\Qp, \alpha^\Z, \beta^\Z)$ is decidable iff $Th(\Qp, \alpha^\Z\beta^\Z)$ is decidable iff Mann property is effective for the group $\alpha^\Z\beta^\Z$ for all $\alpha,\beta\in \N$ with $v(\alpha), v(\beta)$ not both positive. In particular if $\alpha, \beta$ are not multiplicatively independent, then the theory is decidable.
\end{cor}
\begin{proof}
First if $v(\alpha)=v(\beta)=0$ then let us remark that $\alpha^p, \beta^p\in 1+p\Zp$. As $(\Qp, \alpha^\Z, \beta^\Z)$ is definable in $(\Qp, (\alpha^p)^\Z, (\beta^p)^\Z)$ we may therefore assume that $\alpha,\beta\in 1+p\Zp$. Similarly we may assume that $v_p(\alpha-1)=v_p(\beta-1)$. 
\par Let $G=\alpha^\Z$ and $H=\beta^\Z$. If $G\cap H\not=\{1\}$ (iff $\alpha, \beta$ are not multiplicatively independent) then $(\Qp, G,H)$ is definable in $(\Qp, G\cap H)$ (for note that $G\cap H$ has finite index in $G, H$). The theory of this latter structure is decidable by Theorem 2.4 in \cite{Mariaule2}: this theorem axiomatises the theory of $(\Qp, G\cap H)$. All the axioms are obviously recursively enumerable except for the Mann axioms. But as $G\cap H$ is a rank $1$ cyclic group the Mann axioms are effective by \cite{Gunaydin-vdD} Proposition 8.7.
\par Otherwise $G\cap H=\{1\}$ and Theorem \ref{completeness} gives an axiomatisation of $Th(\Qp, G, H)$.	Again it is obvious that all axioms are recursively enumerable except for the Mann axioms. Now we remark that if $Th(\Qp, \alpha^\Z, \beta^\Z)$ is decidable then the collection of Mann axioms for the group $\alpha^\Z \beta^\Z$ is recursively enumerable and conversely. On the other hand by \cite{Mariaule2} Theorem 2.4 it is also the case that $Th(\Qp, \alpha^\Z  \beta^\Z)$ is decidable iff the collection of Mann axioms for the group $\alpha^\Z \beta^\Z$ is recursively enumerable. 
\par Now if $v_p(\alpha)=0$ then $(\Qp, \alpha^\Z, \beta^\Z)$ and $(\Qp, \alpha^\Z\cdot \beta^\Z)$ are bi-interpretable. Indeed, $\alpha^\Z\cdot \beta^\Z\cap \Zp^{\times}=\beta^\Z$. Furthermore the decidability of $Th(\Qp, \alpha^\Z, \beta^\Z)$ is equivalent to effective Mann property for $\alpha^\Z\cdot \beta^\Z$ by Theorem 2.4 in \cite{Mariaule2}.
\end{proof}

\subsection{$Th(\Qp, \alpha^\Z, \beta^\Z)$ is NIP}\label{NIP section}
We will now prove that the theory of $(\Qp, \alpha^\Z, \beta^\Z)$ is NIP for $\alpha, \beta\in \N$ not both with positive valuation. Let us remark that if $v_p(\alpha)=v_p(\beta)=0$ we can assume that $v_p(\alpha-1)=v_p(\beta-1)>0$ like we did in the proof of Corollary \ref{effective MP}. We will tacitly use this reduction in the next results. First we give first three results of quantifier simplification:

\begin{prop}\label{ind structure} Let $(K, G, H)$ be a model of $Th(\Qp, \alpha^\Z, \beta^\Z)$ with $v_p(\alpha)=v_p( \beta)=0$. A subset of $G^m$ if definable iff it is a boolean combination of sets of the forms $X\cap Y$ where $X$ is definable in $(K, GH)$ and $Y\subset G$ is definable in the language of valued groups.
\end{prop}
\begin{proof}
First we prove that $X= X'\cap Y'$ where  $X'$ is definable in $K$ and $Y$ definable in the pair of valued groups $(GH, G)$. Then by quantifier elimination for pairs of valued groups (Theorem \ref{pair of $p$-valued groups}) we can reorganise $X'$ and $Y'$ to obtain the proposition.
\par It is sufficient to prove the following: let $(K_1,G_1,H_1)$ and $(K_2,G_2,H_2)$ be two $|K|^+$-saturated expansion of $(K, G, H)$. Let $\overline{g}_1\in G_1$ and $\overline{g}_2\subset G_2$ such that for any formula $\Psi(\overline{x})$ in the language of rings and parameters in $K$ and for any formula $\varphi(\overline{y})$ in the language of pairs of groups and parameters in $GH$,
$$ (K_1,G_1,H_1)\vDash \Psi(\overline{g}_1)\wedge \varphi(\overline{g}_1) \mbox{ iff } (K_2,G_2,H_2)\vDash \Psi(\overline{g}_2)\wedge \varphi(\overline{g}_2) \qquad (*)$$
Then $tp(\overline{g}_1/(K,G,H))= tp(\overline{g}_2/(K,G,H))$. For it is sufficient to prove that there is an element of the back-and-forth system in the proof of Theorem \ref{completeness} that takes $\overline{g}_1$ to $\overline{g}_2$. As $\overline{g}_1\subset G_1$ we are in case (1) of that proof. That case only use hypothesis $(*)$ to extends the embedding so we are done.
\end{proof}

\begin{definition} Let $T$ be a $\mathcal{L}$-theory, $\mathcal{M}\vDash T$ and $P\subset M$. Let $T_P=Th(M, P)$ in the language $\mathcal{L}\cup\{A\}$ where $A$ is a unary predicate interpreted by $P$. We say that $T_p$ is \emph{bounded} if any formula is equivalent to a boolean combination of formulas of the type
$$\exists \overline{y} (\bigwedge_i P(y_i)\wedge \Phi(\overline{x},\overline{y}),$$
where $\Phi$ is a $\mathcal{L}$-formula (with parameters).
\end{definition}
We show that $Th(\Qp, \alpha^\Z)_{\beta^{Z}}$ is bounded:
\begin{prop}\label{bounded formulas}  Let $(K,G,H)$ be a model of $Th(\Qp, \alpha^\Z, \beta^\Z)$ with $v_p(\alpha)=v_p( \beta)=0$. Every definable subset of $(K,G, H)$ is a boolean combination of subset of $K^n$ defined by formulas $\exists \overline{y}\exists \overline{z} (\overline{y}\subset G \wedge \overline{z}\subset H\wedge \Phi(\overline{x},\overline{y},\overline{z})$ where $\Phi$ is a $\mathcal{L}_{Mac}$-quantifier-free formula.
\end{prop}
\begin{proof} As in the proof of the last proposition it sufficient to prove that for all $(K_1,G_1,H_1)$, $(K_2,G_2,H_2)$ $|K|^+$-saturated expansions of $(K,G,H)$, for all $\overline{x}\in K_1^n$ and $\overline{y}\in K_2^n$ such that $\overline{x}$ and $\overline{y}$ satisfy the same formulas of the type $\exists \overline{y}\exists \overline{z} (\overline{y}\subset G \wedge \overline{z}\subset H\wedge \Phi(\overline{x},\overline{y},\overline{z})$ like in the hypothesis then $tp_{(K, G, H)}(\overline{x})=tp_{(K, G, H)}(\overline{y})$. For it is sufficient to find an embedding $\iota$ in the back-and-forth system in the proof of Theorem \ref{completeness} that takes $\overline{x}$ to $\overline{y}$.
\par Assume that $x_n$ is algebraic over $\Q(G_1H_1)(x_1,\cdots, x_{n-1})$ i.e., there is $g_1,\cdots, g_l\in G_1$, $h_1,\cdots, h_t\in H_1$ and a $\mathcal{L}_{Mac}$formula $\varphi(\overline{u},\overline{v},\overline{w})$ such that
$$(K_1, G_1, H_1)\vDash \varphi(x_n, x_1,\cdots, x_{n-1},\overline{g},\overline{h})\wedge \exists^{\leq n} u_n \varphi(u_n, x_1,\cdots, x_{n-1},\overline{g},\overline{h}).$$
So
$$(K_1, G_1, H_1)\vDash \exists \overline{v}\in G_1 \exists \overline{w}\in H_1 \varphi(x_n, x_1,\cdots, x_{n-1},\overline{v},\overline{w})\wedge \exists^{\leq n} u_n \varphi(u_n, x_1,\cdots, x_{n-1},\overline{v},\overline{w}).$$
Now by assumption
$$(K_2, G_2, H_2)\vDash \exists \overline{v}\in G_2 \exists \overline{w}\in H_2 \varphi(y_n, y_1,\cdots, y_{n-1},\overline{v},\overline{w})\wedge \exists^{\leq n} u_n \varphi(u_n, y_1,\cdots, y_{n-1},\overline{v},\overline{w}).$$
That is there is $g'_1,\cdots, g'_l\in G_2$ and $h'_1,\cdots, h'_t\in H_2$ such that $y_n$ is algebraic over $\Q(y_1,\cdots, y_n, \overline{g'}, \overline{h'})$. By compactness and assumption (extending $\varphi$ if necessary) we may assume that $\overline{g}$ and $\overline{g'}$ (as well as $\overline{h}$ and $\overline{h'}$) satisfies the same formulas of the type $\Psi_1\wedge \Psi_2$ where $\Psi_1$ is a $\mathcal{L}(K, GH)$-formula and $\Psi_2$ is a formula in the language of valued groups. Then as in Proposition \ref{ind structure} we can find an embedding $\iota$ in the back-and-forth system that sends $(\overline{g},\overline{h})$ to $(\overline{g'},\overline{h'})$. So we can assume that $\overline{g},\overline{g'}, \overline{h},\overline{h'} \subset GH$. Now by induction we can assume that $x_1,\cdots, x_r$ are algebraically independent over $\Q(G_1H_1)$. By symmetry also $y_1,\cdots, y_r$ are algebraically independent over $\Q(G_2H_2)$. As $\overline{x}$ and $\overline{y}$ satisfies the same $\mathcal{L}_{Mac}$-formulas we get an isomorphism between $K(x_1,\cdots, x_r)^h$ and $K(y_1,\cdots, y_r)^h$ in the back-and-forth system. As $x_i$ algebraic over $\Q(G, x_1,\cdots, x_r)\subset K(x_1,\cdots, x_r)^h$ for all $i$ (and similarly in $K_2$) we are done.
\end{proof}

\begin{prop}\label{ind structure mixed case} Let $(K, G, H)$ be a model of $Th(\Qp, \alpha^\Z, \beta^\Z)$ with $v_p(\alpha)>0$, $v_p( \beta)=0$.  A subset of $H^m$ if definable iff it is a boolean combination of sets of the forms $X\cap Y$ where $X$ is definable in $(K, G)$ and $Y\subset G$ is definable in the language of valued groups.
\end{prop}
\begin{proof} The proof is similar to Proposition \ref{ind structure}. Here we use back-and-forth system in the proof of Theorem 2.4 in \cite{Mariaule2}.
\end{proof}

\begin{theorem}\label{NIP} $Th(\Qp, \alpha^\Z, \beta^\Z)$ is NIP if $v_p(\alpha)=0$ or $v_p(\beta)=0$.
\end{theorem}
\begin{proof} We will use Corollary 2.5 in \cite{C-S}. For let $T=Th(\Qp, \alpha^\Z)$, $T_P=Th(\Qp, \alpha^Z, \beta^\Z)$ and $h_{ind}=Th(h, (R_\Phi) )$ where $\Phi$ runs over all $\mathcal{L}_{\alpha^\Z}$-formula (with parameters) and $R_\phi$ is a predicate interpreted by $H^n\cap \Phi(K^n)$. Corollary 2.5 in \cite{C-S} states that if $T$ is NIP, $T_P$ is bounded and $H_{ind}$ is NIP then $T_P$ is NIP.
\par First we deal with the case $v_p(\alpha)=v_p(\beta)=0$. By Proposition \ref{bounded formulas} $T_P$ is bounded. By \cite{Mariaule2} Theorem 6.7 $T$ is NIP. It remains to prove that $H_{ind}$ is NIP. For by Proposition \ref{ind structure} it is sufficient to prove that any formula of the type $\Phi\wedge \varphi$ is NIP in $H_{ind}$ where $\Phi$ is a formula in the language of the pair $(\Qp ,\alpha^Z)$ and $\varphi$ is a formula in the language of $p$-valued groups. Let $(a_i; i\in I)$ be an indiscernible sequence in $H_{ind}$ and $b\in H$. Then by definition of the language for this structure $(a_i; i\in I)$ is indiscernible in $(K, G)$. So as $Th(\Qp, \alpha^\Z)$ is NIP, $(\Qp, \alpha^\Z)\vDash \Phi(a_i,b)$ eventually (or $(\Qp, \alpha^\Z)\vDash \neg \Phi(a_i,b)$ eventually). Similarly $(a_i; i\in I)$ is indiscernible in $H$ for the language of valued groups. By \cite{Mariaule2} Theorem 1.7 $Th(H)$ as valued group is NIP. So $H_{ind}\vDash \Phi(a_i,b)\wedge \varphi(a_i,b)$ eventually or $H_{ind}\vDash \neg(\Phi(a_i,b)\wedge \varphi(a_i,b))$ eventually i.e., $H_{ind}$ is NIP.
\par If $v_p(\alpha)>0$ the proof is similar:  $T_P$ is bounded ( \cite{Mariaule2} Proposition 3.3) and $T$ is NIP (\cite{Mariaule2} Corollary 6.5). We use Proposition \ref{ind structure mixed case} as above to prove that $H_{ind}$ is NIP.
\end{proof}

\begin{cor} $Th(\Qp, \alpha^\Z, \beta^\Z)$ does not interpret $(\Z,+,\cdot, 0, 1)$ if $v_p(\alpha)$ or $v_p(\beta)$ is zero.
\end{cor}
\begin{proof} This is immediate from Theorem \ref{NIP} as NIP theories do not interpret Peano arithmetic (in fact any non-NIP theory).
\end{proof}

\bibliographystyle{plain}
\bibliography{Biblio}

\noindent Nathana\"el Mariaule\\
Universit\'e de Mons, Belgium\\
\emph{E-mail address: Nathanael.MARIAULE@umons.ac.be}

\end{document}